\documentclass{amsart}
\usepackage{amsfonts}
\usepackage{graphicx}
\usepackage{amscd}
\usepackage{mathrsfs}
\usepackage{amssymb}
\usepackage{amsmath}
\usepackage{amsthm}

\newtheorem{theorem}{Theorem}[section]
\newtheorem{proposition}{Proposition}[section]

\title[The computation of the M\"{o}bius function of a M\"{o}bius category]{\vspace*{2cm}\it{The computation of the M\"{o}bius function of a M\"{o}bius category}}
\author[Emil Daniel Schwab \& Juan Villarreal] {Emil Daniel Schwab \& Juan Villarreal}
\begin{document}
\topmargin=2.5cm \textwidth=12.5cm \textheight=18cm \maketitle
\begin{center}
{\it Department of Mathematical Sciences\\
The University of Texas at El Paso\\
El Paso, Texas 79968}\\
{\tt eschwab@utep.edu}\\
{\tt jcvillarreal2@miners.utep.edu}
\end{center}
\begin{abstract}
The paper presents some results for reducing the computation of
the M\"{o}bius functon of a M\"{o}bius category that arises from a
combinatorial inverse semigroup to that of locally finite partially ordered sets.
We illustrate the computation of the M\"{o}bius function with an example.


\end{abstract}
\vskip0.5cm
\section{Introduction (M\"{o}bius categories)}
The theory of M\"{o}bius functions in categories was initiated and developed by Leroux and collaborators such as Content and Lemay (see [1],[6],[7]), and in the same time by Haigh [2]. Recently, Leinster [5], Lawvere and Menni [3] have brought attention to the problem of M\"{o}bius inversion in categories in a broader context. In previous papers [9]-[13], the first author of the present paper has found connections between the theory of combinatorial inverse semigroups and the theory of M\"{o}bius categories. The combinatorial inverse semigroups provide special examples of M\"{o}bius categories.

A M\"{o}bius category $C$ (in the sense of Leroux) is a small category satisfying the following conditions:

(1) any morphism $f\in{MorC}$ has only a finitely many non-trivial factorizations $f=gh$;

(2) an incidence function $\xi:MorC\rightarrow{\mathbb{C}}$ has a convolution inverse if and only if $\xi(f)\neq0$ for each identity morphism $f$ of $C$. (The convolution $\ast$ of two incidence functions $\xi$ and $\eta$ is defined by: $(\xi\ast\eta)(f)=\sum_{f=gh}\xi(g)\eta(h)$, and the convolution identity is $\delta$ given by: $\delta(f)=1$ if $f$ is an identity morphism and $\delta(f)=0$ otherwise.)

The M\"{o}bius function $\mu$ of a M\"{o}bius category $C$ is the convolution inverse of the zeta function $\zeta$: $\zeta(f)=1$ for any morphism $f$ of $C$. The M\"{o}bius inversion formula is then nothing but the statement: $\xi=\eta\ast\zeta\Leftrightarrow\eta=\xi\ast\mu$. This is also the M\"{o}bius inversion formula in number theory, the functions being arithmetic functions and the convolution $\ast$, the Dirichlet product.

Now, a partially ordered set (poset) $(P,\leq)$ as a category (the objects are the elements of $P$, and there is a morphism $x\rightarrow{y}$ if and only if $x\leq{y}$) is M\"{o}bius if and only if it is locally finite (i.e. any interval in $P$ is finite). So, M\"{o}bius inversion, a useful tool in number theory, was generalized to categories via Rota's theory of M\"{o}bius functions. The computation of poset's M\"{o}bius function has a central place in Rota's theory ([8]). There are many remarkable fruitful methods for computing poset's M\"{o}bius functions.

Our main interest here is in the M\"{o}bius function for M\"{o}bius categories that arise from combinatorial inverse semigroups. There are some approaches to compute such M\"{o}bius functions using poset's M\"{o}bius functions. In Section 3 an example is considered to illustrate one of the presented methods.

\vskip0.5cm
\section{Links with inverse semigroups and poset's M\"{o}bius functions}
A semigroup $S$ is an inverse semigroup if every element $s\in{S}$ has a unique inverse $s^{-1}$, in the sense that $ss^{-1}s=s$ and $s^{-1}ss^{-1}=s^{-1}$. An inverse semigroup $S$ possesses a natural partial order relation $\leq$ defined by: $s\leq{t}\Leftrightarrow{s=ss^{-1}t}$. An inverse semigroup is locally finite if the poset of idempotents $(E(S),\leq)$ is locally finite. Two elements $s,t\in{S}$ are said to be $\mathscr{D}$-equivalent if and only if there exists an element $x\in{S}$ such that $s^{-1}s=x^{-1}x$ and $xx^{-1}=tt^{-1}$. A combinatorial semigroup is a semigroup whose all subgroups are trivial.

The natural partial order on inverse semigroups implies a combinatorial approach of inverse semigroups via Rota's M\"{o}bius inversion. The idea was explored successfully by Steinberg [14],[15].

The theory of Leech's ([4]) division categories open a new way to combinatorial approach of inverse semigroups, in this time, via categorical M\"{o}bius inversion. By [9, Theorem 3.3],
\begin{quote}
"the reduced division category $C_{F}(S)$ relative to an idempotent transversal $F$ of the $\mathscr{D}$-classes of an inverse monoid $S$ with $1\in{F}$ is a M\"{o}bius category if and only if $S$ is locally finite and combinatorial".
\end{quote}
This category $C_{F}(S)$ is defined by:
\begin{itemize}
\item $ObC_{F}(S)=F$;
\item $Hom_{C_{F}(S)}(e,f)=\{(s,e)|s\in{S},\ s^{-1}s\leq{e}\ and\ ss^{-1}=f\}$;
\item The composition of two morphisms $(s,e):e\rightarrow{f}$ and $(t,f):f\rightarrow{g}$ is given by $(t,f)\cdot(s,e)=(ts,e)$.
\end{itemize}
The  M\"{o}bius function $\mu$ of this  M\"{o}bius category $C_{F}(S)$ is called the  M\"{o}bius function of $S$. If $S$ is with zero or has no identity element then a slight correction is needed to define the  M\"{o}bius category of $S$ (see [11]).

Starting with a locally finite combinatorial inverse monoid $S$, we are considering three rules for determining $\mu(s,e)$:

(1) The collection of quotients in $C_{F}(S)$ of an object $e$ forms a locally finite partially ordered set in a natural way. The first rule uses the  M\"{o}bius function of this locally finite partially ordered set. For any $e\in{ObC_{F}(S)}$, the set $Q(e)$ of all quotient objects of $e$ with tha canonical partial order $\preceq$ is a locally finite lattice, and (see [9, Theorem 3.5 ($i$)):
$$\mu(s,e)=\mu_{Q(e)}((s,e),(e,e)),$$
for any morphism $(s,e)$ of $C_{F}(S)$, where $\mu_{Q(e)}$ is the  M\"{o}bius function of $(Q(e),\preceq)$.

(2) The natural partial order $\leq$ on $S$, restricted to idempotents is given by: $e,f\in{E(S)}:\ e\leq{f}\Leftrightarrow{e=ef=fe}.$ The second rule for determining the value of the function $\mu$ at $(s,e)$ uses th  M\"{o}bius function $\mu_{E(eSe)}$ of the locally finite lattice $(E(eSe,\leq)$. We have:
$$\mu(s,e)=\mu_{E(eSe)}(s^{-1}s),$$
for any morphism $(s,e)$ of $C_{F}(S)$ (see [9, Theorem 3.5 ($ii$)).

(3) The third rule uses Lawvere intervals. If $f$ is a morphism of a small category $C$ then the Lawvere interval $I(f)$ of $f$ is a category with the set of factorizations of $f$ as objects. The morphisms of $I(f)$ are morphism of $C$ which are compatible with factorizations of $f$. That is, if $f=uv$ and $f=u'v'$ in $C$ then $h:Domg(=Codomh)\rightarrow{Domg'(=Codomh')}$ is a morphism in $I(f)$ from $uv$ to $u'v'$ if $uh=u'$ and $hv'=v$. The composition in $I(f)$ is the same as in $C$.

\begin{theorem}
([2]) A small category $C$ is M\"{o}bius if and only if all Lawvere intervals of $C$ are finite and one-way (i.e. $Hom_{C}(X,Y)\neq\emptyset$, $Hom_{C}(Y,X)\neq\emptyset\Rightarrow{X=Y}$; and $|Hom_{C}(X,X)|=1$ for any object $X$ in $C$).
\end{theorem}

\begin{theorem}
([2],[12]) If $C$ is M\"{o}bius then any Lawvere interval is M\"{o}bius and for any morphism $f:X\rightarrow{Y}$ in $C$,
$$\mu(f)=\mu_{I(f)}(f),$$
where $\mu_{I(f)}$ is the M\"{o}bius function of $I(f)$. (The morphism $f$ of $C$ is also a morphism of $I(f)$ from $f1_{X}$ to $1_{Y}f$; and if $I(f)$ is a bounded poset then $\mu(f)=\mu_{I(f)}(0,1)$ with $0=f1_{X}$ and $1=1_{Y}f$.)
\end{theorem}

\begin{theorem}
([12]) If a small category $C$ is M\"{o}bius that arises from a combinatorial inverse semigroup $S$ then any Lawvere interval is a finite lattice (as a category), and for any morphism $f$,
$$\mu(f)=\mu_{I(f)}(0,1),$$
where $\mu_{I(f)}$ is the M\"{o}bius function of the finite lattice $I(f)$, $0$ is the least element and $1$ is the greatest element of $I(f)$.
\end{theorem}

In the next section we illustrate the computation of the M\"{o}bius function with an example. We consider a category $C_{m}$ and based on the above theorems:

1) we show that $C_{m}$ is M\"{o}bius;

2) we find the M\"{o}bius funtion of $C_{m}$.

We focus our attention on the above problems and we will look at the starting inverse semigroup only at the end of the paper.
\vskip0.5cm
\section{An example}
Let $m$ be a positive integer ($m>1$), $Z_{m}$ the cyclic group of addition modulo $m$, $Z_{-}$ the set of non-positive integers, and $Z_{+}$ the set of non-negative integers. Now, let $C_{m}$ be the category defined by:

\begin{itemize}
\item $ObC_{m}=Z_{m}\times{Z_{-}}$
\item $Hom_{C_{m}}((\overline{x},i),(\overline{y},j))=\{(a,\overline{x},i,j)|a\in{Z_{+}},\ a\leq{i-j},\ \overline{a}+\overline{x}=\overline{y}\}$
\item $(b,\overline{y},j,k)\circ(a,\overline{x},i,j)=(a+b,\overline{x},i,k)$ is the composition of two morphisms $(a,\overline{x},i,j):(\overline{x},i)\rightarrow{(\overline{y},j)}$ and $(b,\overline{y},j,k):(\overline{y},j)\rightarrow{(\overline{z},k)}$.
\end{itemize}

The first question is this: is the above category M\"{o}bius ? We will answer this question by applying Theorem 2.1. The second question concerns the computation of the M\"{o}bius function $\mu$ of $C_{m}$.

Let $(a,\overline{x},i,j):(\overline{x},i)\rightarrow{(\overline{y},j)}$ be a morphism in $C_{m}$. First, we examine the factorizations in $C_{m}$ of this morphism, that is the objects of the category (of the Lawvere interval) $I(a,\overline{x},i,j)$. Let
$$(a,\overline{x},i,j)=(a-b,\overline{z},k,j)\circ{(b,\overline{x},i,k)}$$
be a factorization of $(a,\overline{x},i,j)$, i.e. the Diagram 1 is commutative.


\begin{figure}[!h]\centering
\includegraphics*[width=9cm]{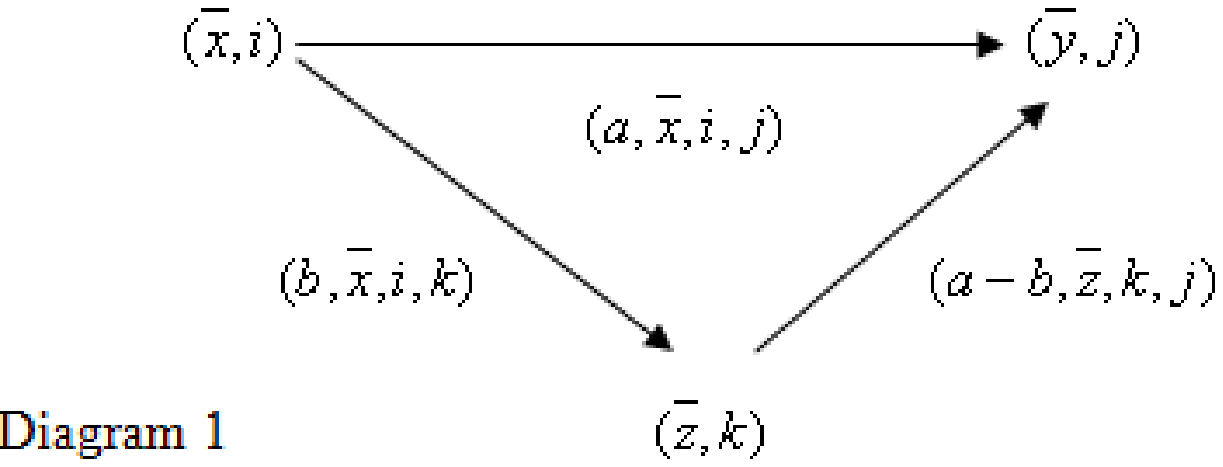}
\end{figure}

We have:
$$a>0,\ i,j\leq0,\ a\leq{i-j}\ and\ \overline{a}+\overline{x}=\overline{y}\  (since\  (a,\overline{x},i,j)\in{Hom_{C_{m}}(\overline{x},i),(\overline{y},j)})$$
and
$$0\leq{b}\leq{a},\ k\leq0,\ b\leq{i-k},\ a-b\leq{k-j},\ \overline{b}+\overline{x}=\overline{z}\ (\overline{a-b}+\overline{z}=\overline{y}\ is\ a\ consequence).$$

Now, for a fixed integer $b$, $0\leq{b}\leq{a}$,

\begin{itemize}
\item $\overline{z}$ is uniquely determined by: $\overline{z}=\overline{b}+\overline{x}$. We denote by $\overline{z_{b}}$ this residue class.
\item The values of $k$ are determined by the condition: $a-b+j\leq{k}\leq{i-b}.$ Since,
$$i-b-(a-b+j)=i-j-a\geq0,$$
it follows that the values of $k$ are the following:
$$k_{0}=a-b+j,\ k_{1}=a-b+j+1,...,k_{t}=a-b+j+t,...$$
$$...,k_{i-j-a}=i-b.$$
\end{itemize}

Thus,

\begin{proposition}
The set of objects of $I(a,\overline{x},i,j)$ is finite for any morphism $(a,\overline{x},i,j)$ of $C_{m}$.
\end{proposition}

With the above notations there exists a morphism of $C_{m}$ from $(\overline{z_{b}},k_{t})$ to $(\overline{z_{b}},k_{p})$ such that the Diagram 2 is commutative if and only if $k_{p}\leq{k_{t}}$ (that is, if and only if $p\leq{t}$), where $0\leq{p,t}\leq{i-j-a}$. It is clear that this morphism $(0,\overline{z_{b}},k_{t},k_{p})$ of $C_{m}$ is uniquely determined by the commutative Diagram 2.

\begin{figure}[!h]\centering
 \includegraphics*[width=9cm]{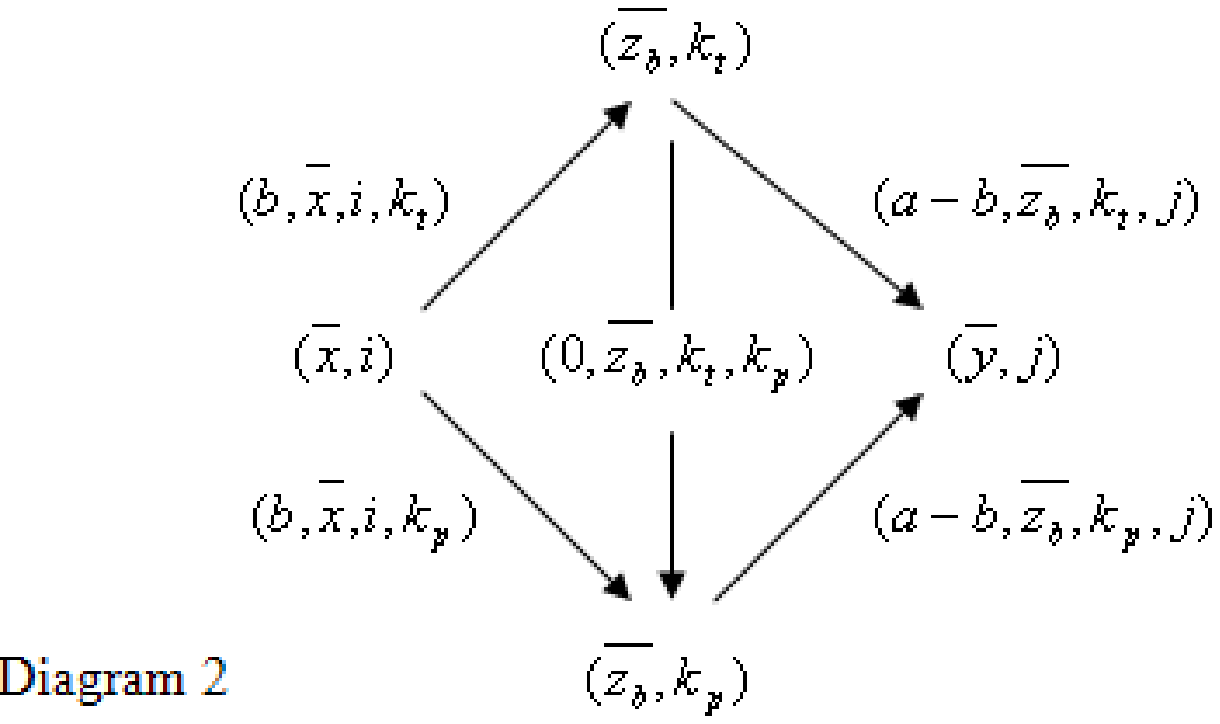}
\end{figure}

Thus,

\begin{proposition}
The set of morphisms $Hom_{I(a,\overline{x},i,j)}(X_{b,t},X_{b,p})$, where
\[
\stackrel{\ \ \ \ \ \ (b,\overline{x},i,k_{t})}{X_{b,t}=(\overline{x},i)\ \ \longrightarrow\ \ \ \ (\overline{z_{b}},k_{t})}
\stackrel{(a-b,\overline{z_{b}},k_{t},j)\ \ \ \ }{\ \ \longrightarrow{\ \ \ \ (\overline{y},j)}}
\]
and
\[
\stackrel{\ \ \ \ \ \ (b,\overline{x},i,k_{p})}{X_{b,p}=(\overline{x},i)\ \ \longrightarrow\ \ \ \ (\overline{z_{b}},k_{p})}
\stackrel{(a-b,\overline{z_{b}},k_{p},j)\ \ \ \ }{\ \ \longrightarrow{\ \ \ \ (\overline{y},j)}}
\]
is non-empty and it is a singleton if and only if $p\leq{t}$.
\end{proposition}


\begin{figure*}
\begin{center}
\includegraphics[width=9cm]{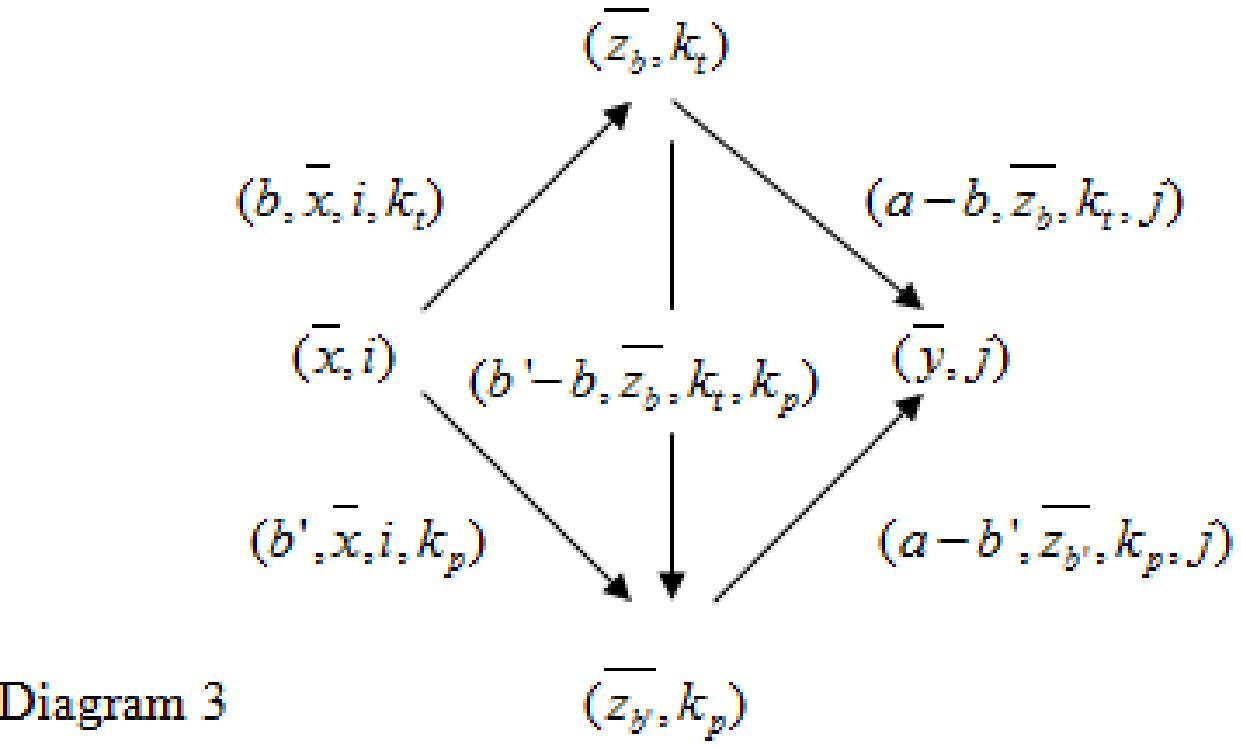}
\end{center}
\end{figure*}

A similar examination of the commutativity in $C_{m}$ of the Diagram 3 leads us to conclude that (the notations are the same as in Proposition 3.2):

\begin{proposition}
The set of morphisms $Hom_{I(a,\overline{x},i,j)}(X_{b,t},X_{b',p})$ is non-empty and it is a singleton if and only if $b\leq{b'}$ and $p\leq{t}$.
\end{proposition}
By Theorem 2.1 and the above results, we obtain the following proposition:
\begin{proposition}
The category $C_{m}$ is a M\"{o}bius category, and every Lawvere interval of $C_{m}$ is a finite lattice.
\end{proposition}

The Diagram 4 is the Hasse diagram of the lattice $I(a,\overline{x},i,j)$.

\begin{figure*}
\begin{center}
\includegraphics[width=15cm]{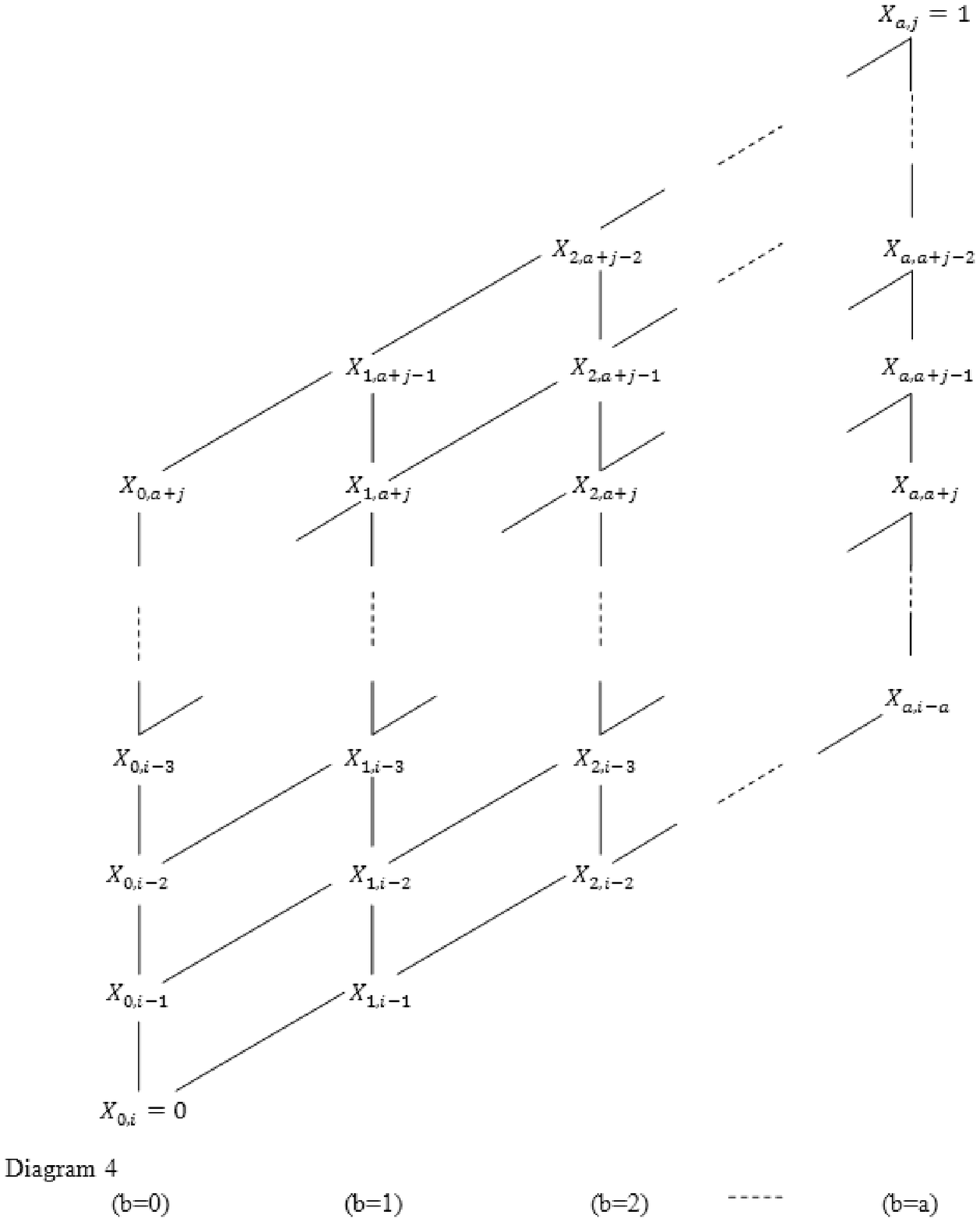}
\end{center}
\end{figure*}


Now, a routine evaluation on the Hasse Diagram 4 of the poset's M\"{o}bius function $\mu_{I(a,\overline{x},i,j)}$, and the equalities $\mu(a,\overline{x},i,j)=\mu_{I(a,\overline{x},i,j)}(a,\overline{x},i,j)=\mu_{I(a,\overline{x},i,j)}(0,1)$ implies:

\begin{proposition}
The M\"{o}bius function $\mu$ of the M\"{o}bius category $C_{m}$ is given by:
  \[\mu(a,\overline{x},i,j)=\left\{\begin{array}{rcl}
1&\mbox{if}&a=0\ and\ j=i\ \ \ \ \ \ \ \ or\ \ \ a=1\ and\ j=i-2\\
-1&\mbox{if}&a=0\ and\ j=i-1\ \ \ or\ \ \ a=1\ and\ j=i-1\\
0&\mbox{\ }&\ \ \ \ \ \ \ \ \ \ \ \ \ \ \ \ \ \ \ \ \ \ \ otherwise.
\end{array}\right. \]
\end{proposition}

\vskip0.5cm
\section{Final remarks}
The Theorem 2.3 suggests that the M\"{o}bius category $C_{m}$ arises from a combinatorial inverse semigroup. An examination of the M\"{o}bius function leads us to the free monogenic inverse monoid (see [10] , Proposition 3.3]). But $C_{m}$ is not the M\"{o}bius category (i.e. the reduced division category) of the free monogenic inverse monoid. The category $C_{m}$ is the M\"{o}bius category of a $\varrho$-semigroup of the free monogenic inverse monoid, namely of the one-dimensional tiling semigroup in the periodic case: the period has length $m$ and involves each tile exactly once (see [13]). By a change from a combinatorial inverse monoid to his $\varrho$-semigroup, the M\"{o}bius function becomes an invariant ([13, Theorem 2.4]).

Now, another way to compute the M\"{o}bius function of a M\"{o}bius category is given below (using a comparison method). Let $D_{m}$ be the M\"{o}bius category defined by:

\begin{itemize}
\item $ObD_{m}=Z_{m}$
\item $Hom_{D_{m}}(\overline{x},\overline{y})=\{(\alpha,\overline{x})\in{Z_{+}\times{Z_{m}}}|\alpha\geq{x}$ and $\alpha\equiv{y}$ (mod$m$)\};
\item If $(\alpha,\overline{x})\in{Hom_{D_{m}}}(\overline{x},\overline{y})$ and $(\beta,\overline{y})\in{Hom_{D_{m}}}(\overline{y},\overline{z})$ then  $(\beta,\overline{y})\cdot{ (\alpha,\overline{x})}=(\beta-y+\alpha,\overline{x})$.
\end{itemize}
The functor $F:C_{m}\rightarrow{D_{m}}$ given by:
 \begin{itemize}
\item $F(\overline{x},i)=\overline{x}$
\item $F(a,\overline{x},i,j)=(a+x,\overline{x})$
\end{itemize}
is surjective that commute with Lawvere intervals $I$. The chain of Diagram 5 is the image through the functor $F$ of the Lawvere interval $I(a,\overline{x},i,j)$ (Diagram 4).



\begin{figure}[!h]\centering
 \includegraphics*[width=11cm]{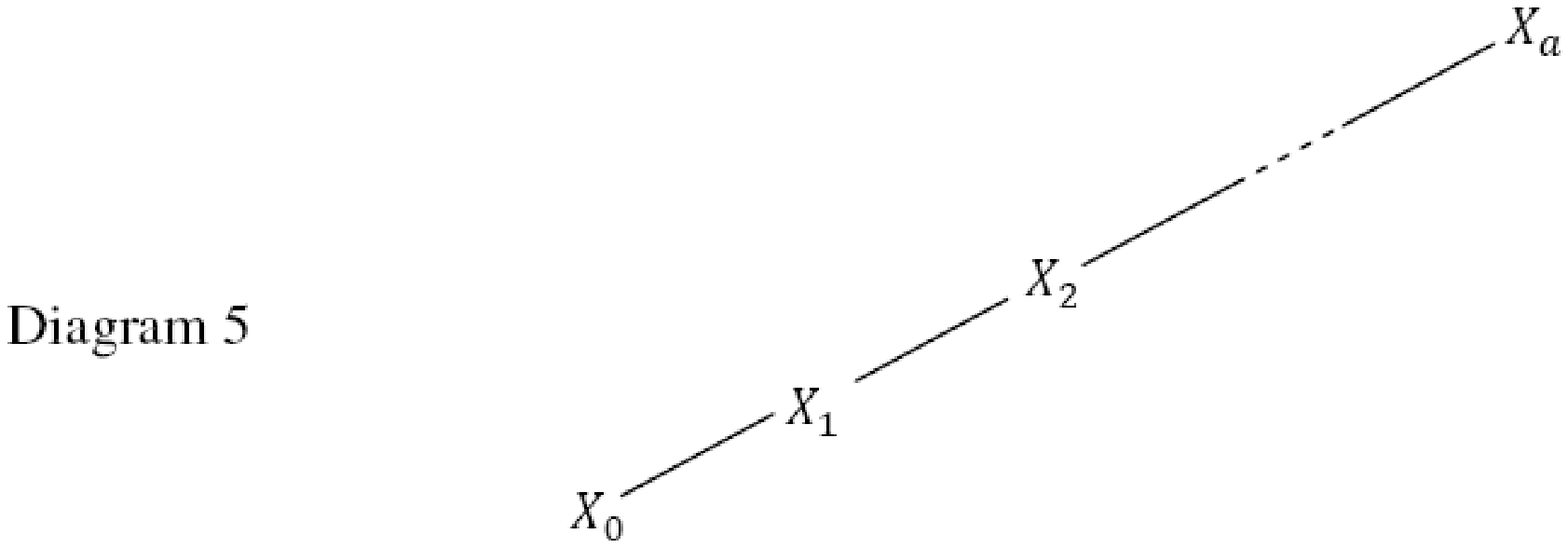}
\end{figure}

It follows that (where $\alpha=a+x$):
\[\mu(\alpha,\overline{x})=\left\{\begin{array}{rcl}
1&\mbox{if}&\alpha=x\\
-1&\mbox{if}&\alpha=x+1\\
0&\mbox{if}&\alpha\geq{x+2}
\end{array}\right. \]
$\mu$ being the M\"{o}bius function of the M\"{o}bius category $D_{m}$.

\vskip0.5cm

\end{document}